\documentclass[12pt]{article}
\newcommand{\average}{-\!\!\!\!\!\!\int}

\newcommand{\comment}[1]{}

\renewcommand{\div}{{\mathop{\rm\,div\,}\nolimits}}

\newcommand{\dist}{\mathop{\rm dist}\nolimits}

\providecommand{\qed}{\vrule height 6pt depth 0pt width 3 pt}

\newcommand{\reals}{{\bf R}}

\newcommand{\supp}{\mathop{\rm supp}\nolimits }

\newcommand{\BibTeX}{{\rm B\kern-.05em{\sc i\kern-.025em b}\kern-.08em     
    T\kern-.1667em\lower.7ex\hbox{E}\kern-.125emX}}

\usepackage{hyperref}

\newcommand{\note}[1]{}
\renewcommand\marginpar[1]{}

\usepackage{amsmath}
\usepackage{amsmath} 
\usepackage{amssymb}

\newcommand{\locdom}[2]{{\Omega_{#2}(#1)}}
\newcommand{\sball}[2]{{\Delta_{#2}(#1)}}

\newcommand{\avg}[3]{ \bar #1 _{#2,#3}}

\newenvironment{proof}[1][Proof]{\begin{trivlist}\item[\hskip \labelsep
{\it #1. }]}{\hfill \qed \goodbreak \end{trivlist}}   
\newenvironment{remark}[1][Remark]{\begin{trivlist}\item[\hskip \labelsep
{\it #1. }]}{  \goodbreak \end{trivlist}}   
\numberwithin{equation}{section} 
\newtheorem{theorem}[equation]{Theorem}
\newtheorem{proposition}[equation]{Proposition}
\newtheorem{corollary}[equation]{Corollary}

\renewcommand{\theequation}{\arabic{section}.\arabic{equation}}

\usepackage{amssymb}

\begin{document}
\title{The Green function for the mixed problem for the linear
  Stokes system in domains in the plane} 

\author{
K.A.~Ott\footnote{Katharine Ott is partially supported by a grant from the
  U.S.~National Science Foundation, 
  DMS 1201104.} 
\\Department of Mathematics\\University of Kentucky\\Lexington,
KY 40506-0027, USA
\and 
 S. Kim
\footnote{
Seick Kim is supported by NRF Grant No. 2012-040411.}
  \\ Department of Mathematics \\ Yonsei University \\ Seoul, 120-749, Korea
\and
 R.M.~Brown\footnote{
Russell Brown   is  partially supported by a grant from the Simons
Foundation (\#195075).
} \\ Department of Mathematics\\ University of Kentucky \\
Lexington, KY 40506-0027, USA}

\date{}

\maketitle

\abstract{ 
We construct the Green function for the mixed boundary value problem
for the  linear  Stokes system in a
two-dimensional Lipschitz domain. 
  }

\section{Introduction}\label{Introduction}

 Let $ \Omega \subset \reals ^2$ be a 
domain and suppose that we have a decomposition of the boundary $
\partial \Omega = D\cup N$ with $ D  \cap N = \emptyset$. 
We consider the mixed  boundary value problem for the linear Stokes system
\begin{equation}\label{MP}
\left\{ \begin{aligned} 
&-\Delta u + \nabla p = f , \qquad &&\mbox{in } \Omega \\
& - \div u = g , \qquad            &&\mbox{in } \Omega \\
&2 \nu \epsilon(u)  - p \nu = f _N,\qquad  &&\mbox{on } N \\
&u = f _D , \qquad                 &&\mbox{on }  D .
\end{aligned}
\right.
\end{equation}
In the boundary value problem (\ref{MP}),  the functions $f$, $f_D$,
$f_N$ and $g$ are given 
and we look for a vector-valued function $ u :  \Omega \rightarrow
\reals^2$ and scalar function $p : \Omega \rightarrow \reals$ which
satisfy the above conditions.   We use $ \epsilon(u) $ to denote the
symmetric part of the gradient of $u$ and $ \nu$ is the outer unit
normal to $ \partial \Omega$. 

The goal of this note is to  give conditions on the decomposition
$\partial \Omega = D\cup N$ 
that allow us to  construct the Green function for the
boundary value problem (\ref{MP})
in a Lipschitz domain in two dimensions. Our argument
begins with an idea of D.~Mitrea and I.~Mitrea \cite{MR2763343} who
construct Green functions for second order elliptic equations with
Dirichlet boundary conditions in two dimensions by extrapolating from
the standard theory of weak solutions with gradient in $L^2$.  We
expect the Green function to have a gradient in the Lorentz space $L^
{ 2, \infty }$, or weak $L^2$, and an argument that involves
perturbing the function space allows us to
extend the existence theory from $L^2$ to weak $L^2$. 

Once we have the existence of a solution with gradient in $L^ { 2,
\infty}( \Omega)$, we use local regularity estimates for solutions
to establish pointwise estimates for the Green function. These
estimates are  mainly of  interest near the boundary since 
 the interior regularity of solutions is
well-understood.

The approach we use is  limited to two dimensions. To study the Green
function for a  
boundary value problem for a second order operator  
in an $n$-dimensional domain, one would need to
study solutions that have their gradient in $L^ { n / ( n -1 ) ,
\infty}$. When $ n =2$, $ n /(n-1) =2$ and it is not
difficult to pass from the standard theory of weak solutions to results for functions that have a gradient  in $L^{2, \infty}$.
When $ n>2$, the gap from $L^2$ to $L^ { n/(n-1) , \infty}$ is too large for this approach to be fruitful.  
  There is a well-known path  to construct the
Green function for the Stokes operator or elliptic systems with
Neumann or Dirichlet boundary conditions in Lipschitz domains in three
dimensions.  This begins with work of Pipher and Verchota
\cite{MR94j:35030} and Dahlberg and Kenig \cite{DK:1990} and continues
in many papers. The recent of work of D.~Mitrea and I.~Mitrea
\cite{MR2763343} includes a construction of the Green function with
Dirichlet boundary conditions for elliptic systems in three
dimensions.  However, the range of exponents for which we can study
the mixed boundary value problem in the non-tangential sense is
smaller than the range for other boundary value problems  and thus we are not able to adapt these arguments to construct the Green function for the mixed problem in three dimensions. 
Furthermore, we are interested in understanding the Green function as a
step towards studying the  mixed problem in the non-tangential sense. 

We recall only a few high points in the study of Green functions for
elliptic operators. 
A classic result of Littmann, Stampachhia and Weinberger
\cite{LSW:1963} gives pointwise estimates for the fundamental solution
(or Green function) 
of an elliptic equation in the plane with bounded and measurable
coefficients. 
 Dolzmann and M\"uller
\cite{MR1354111} construct the Green function in a $C^1$-domain for an
elliptic operator with continuous coefficients. 
Auscher and collaborators \cite{MR1600066} consider elliptic
equations with complex coefficients. Such operators may also be viewed
as elliptic systems. 
Dong and Kim
\cite{MR2485428} construct the Green function for elliptic systems for
operators with bounded and measurable coefficients and with
Dirichlet boundary conditions. 

There is a large body of work related to the study of the Stokes
equations in Lipschitz  domains.  Fabes, Kenig, and Verchota
\cite{FKV:1988} treat the Neumann and Dirichlet problems. They
establish the existence of  solutions with 
the non-tangential maximal function of the
gradient in $L^2( \partial \Omega)$. The results for the
Dirichlet problem were extended to an optimal range of $L^p$-spaces by
Z.~Shen \cite{MR1223521}.  The work of D.~Mitrea and I.~Mitrea \cite{MR2763343}
mentioned above gives a construction of the Green function for the
Stokes operator with Dirichlet boundary conditions in   two and three
dimensions.  The  
mixed problem for the Laplacian  
in Lipschitz domains was the subject of a problem posed
by Kenig in his CBMS lectures \cite{CK:1994}. Recent progress may be
found in the article of Taylor, Ott, and Brown \cite{MR3034453}. This
work and related results on the Lam\'e system \cite{MR3040944} in two
dimensions rely on estimates for the Green function with mixed
boundary conditions. One motivation  for the present note  is
an interest in 
developing the properties of the Green function that are needed to attack
the mixed problem for the linear Stokes system in two dimensions.

There is also substantial interest in studying the mixed problem for
the linearized Stokes equations in
polyhedral domains and obtaining optimal regularity results. Many
polyhedral domains are also Lipschitz domains, however the class of
polyhedral domains includes domains that are not Lipschitz in the
sense defined below, at least in dimension three and higher.  Lipschitz
domains, of course, include many domains that are not polyhedral and
are of interest because the class of Lipschitz domains includes domains
with interesting features at all length scales.  Our treatment of the
mixed problem includes conditions on the decomposition of the boundary
that are scale invariant as well.  We refer the reader to Maz'ya and
Rossmann \cite{MR:2007,MR2563641} and the references cited therein for
additional background on the mixed problem for the Stokes system  in polyhedral
domains.

Finally, we note that there has been recent work on the mixed problem
in domains that are more general than Lipschitz. See work of Auscher,
Badr, Haller-Dintelmann, and Rehberg \cite{ABHR:2012}, 
Haller-Dintelmann, Jonsson, Knees, and Rehberg \cite{HJKR:2012}, and
Brewster, D.~Mitrea, I.~Mitrea, and M.~Mitrea \cite{BMMM:2012}.
Roughly speaking, one needs the set where Neumann data is specified to be
Lipschitz, while the Dirichlet set is allowed to be more general. The
work of Brewster {\em et. al. }has a weaker condition near $N$. It
would be interesting to construct Green functions for the mixed problem
in a similar setting.  One impediment to carrying out the work
reported here in a more general setting is the difficulty of treating
the equation $ \div u = f$ as in Proposition \ref{Bogo}. 

In section \ref{WeakSolutions} we give a  standard weak formulation of 
the mixed problem in (\ref{WMP}).    Using this notion of a  weak
solution 
we
 define the {\em Green function }for the boundary value problem
(\ref{MP}). The Green function is a pair $ (G(x,y), \Pi(x,y))$ where $
G^ { \alpha \beta }: \Omega \times \Omega \rightarrow \reals $ and $
\Pi ^ \alpha : \Omega \times \Omega \rightarrow \reals$ with $ \alpha,
\beta =1,2$. If $(u,p)$ is a weak solution (as defined below)  of the mixed problem with data 
$f$ and $g$ taken from $ C^ \infty _ c ( \Omega)$ and the boundary data $
f_D$ and $f_N$ are zero, then the solution $u$ is given by
\[
u ^ \alpha (x) = \int _ \Omega G^{ \alpha \beta } ( x,y ) f^ \beta (y)
 + \Pi ^ \alpha (x,y ) g(y ) \, dy . 
\]
Since we have uniqueness of weak solutions, it is immediate that the
Green function is unique. 
\note{ Give up on representation formula for the pressure. Leave this
  as an open problem for now. }

Our main result is the following theorem.  The reader will need to
refer to section \ref{Defs} for a  detailed statement of our conditions
on the domain $ \Omega$ and the decomposition of the boundary $
\partial \Omega = D \cup N$.   In the theorem below, we use $d$ to
denote the diameter of $ \Omega$. 

\begin{theorem} 
\label{Green}
Let $\Omega $ be a Lipschitz domain and suppose that $ D$ and $ N$
satisfy the conditions (\ref{AhlDav}), (\ref{DOpen}), and
(\ref{NOpen}).
 There exists a Green function $(G, \Pi)$ for the
boundary value problem (\ref{MP}) and the Green function satisfies the
following:

\begin{align}
& |G(x,y) | \leq C ( 1 + \log ( d / |x-y |)   )
\label{GreenUp} \\
& \| \nabla G(x, \cdot ) \|_ { L^ { 2, \infty } ( \Omega ) } + \|
\Pi (x, \cdot ) \|_ { L^ { 2, \infty } ( \Omega)} \leq C   \label{GreenLorentz} \\
& \nabla_y G(x, \cdot ) , \Pi(x, \cdot ) \in  L^ q ( \Omega
\setminus \Omega _ \rho( x))  \text{ for  } \rho > 0, \frac 1 q> \frac
1 2 -\kappa \label{GreenGrad} \\
& |G(x,y ) - G(x, z ) | \leq  C \left ( \frac { |y -z | } { |x-y | }
\right) ^ \gamma , \qquad |x-y| > 2 |  y-z | \label{GreenHolder} \\
& G^ { \alpha \beta } (x,y) = G^ { \beta \alpha } (y,x) . \label{GreenSymm}
\end{align}
The parameter $ \kappa>0 $ depends on the Lipschitz character of
$\Omega$ and  the H\"older exponent $ \gamma$ depends only on $M$. 
\note{define $d$}

Furthermore  if $u$ is the weak solution of (\ref{WMP}) with  $ f_N \in L^ t ( N)$, $f \in L^ t (
\Omega )$, and $g \in L^ {t} _{1, \partial \Omega}( \Omega)$ for some $t>1$, then we have
\begin{multline}
\label{GreenRep}
u^ \alpha (x) = 
 \int _ \Omega G^ { \alpha \beta } (x,y) ( f^ \beta (y)
- \frac { \partial g } {\partial y _ \beta } (y) )
+ \Pi ^ \alpha(x,y)  g(y)  \, dy 
\\
 + \int _N G^ { \alpha \beta }(x,y)  f_N ^ \beta
(y) \, d\sigma. 
\end{multline}
\end{theorem}

The paper will proceed in the following manner. 
In section \ref{Defs} we introduce the function spaces and machinery
needed to construct the Green function. 
Section \ref{WeakSolutions} establishes the existence of weak
solutions. Section \ref{Local} gives the local regularity needed to
establish the pointwise estimates for the Green function and section
\ref{TheEnd} provides the details of the proof of Theorem \ref{Green}. 

\section{Notations and definitions}
\label{Defs}
\subsection{Domains}

We assume that $ \Omega \subset \reals^2 $ is a {\em Lipschitz domain}. 
Thus $ \Omega$ is a bounded, connected, open set and if $ x\in
\partial\Omega$, then the boundary near $x$ is given by
the graph of a Lipschitz function. More precisely, this means that we
have  constants $M >0$ and $ R_0>0$ so that for each $ x \in \partial
\Omega$, we may find a Lipschitz function $ \phi : \reals \rightarrow
\reals$  such that 
\begin{align}
\label{LipDef1}
\Omega \cap Z_{200 R_0 } (x) & = \{ y : y _ 2 > \phi (y_1) \}  \cap Z _
       {200 R _ 0} (x) \\
\label{LipDef2}
\partial \Omega \cap Z_{200 R_0 } (x) & = \{ y : y _ 2 = \phi (y_1) \}  \cap Z _
       {200 R _ 0} (x),
\end{align}
where for $x=(x_1,x_2) \in \reals^2$,  $ Z _ \rho ( x) = \{ y : |x_1 - y _1 | < \rho,
|x_2- y _ 2 | < (4M+2) \rho\} $ is a {\em coordinate cylinder } centered at $ x\in \partial\Omega $. 
The coordinate system in (\ref{LipDef1}-\ref{LipDef2}) is assumed to
be a rotation of the standard coordinate system. 
Since the domain is bounded, we may  fix a covering of $ \partial \Omega$
by a finite collection of coordinate cylinders $ \{ Z_{R_0}(x_1),
\dots,  Z_{R_0}(x_N)\}$ and we will use these cylinders in the
constructions  below. 

Next we define {\em boundary intervals }$\sball x \rho \subset  \partial \Omega$. 
If $ 0 < \rho < 100R_0$ and $ x \in 
\partial \Omega$, then we set $ \sball x \rho  = Z_ \rho (x) \cap
\partial \Omega$. We also define {\em local domains } $ \Omega_\rho (x) \subset \Omega $.
 These sets will be disks in the interior of
$ \Omega$ and will provide a convenient family of sets for studying
the local regularity of solutions of a boundary value problem near the
boundary. For $ 0 < \rho < 100R_0$, if the distance from $ x $ to $
\partial \Omega$, $ \dist(x, \partial \Omega) > \rho$, then we let $
\locdom x \rho = \{ y : |x-y | < \rho \}$ be the disk centered at $x$
and of radius $ \rho$. If $ \dist(x, \partial \Omega) \leq 
 \rho$, then $x$ lies in some coordinate cylinder $ Z$  and we write 
$ x = \hat x + s e_2$ with $ \hat x \in \partial \Omega\cap Z$   
in the coordinate system of $Z$ and $s>0$. 
We define the local domain  $ \locdom x  \rho = \Omega \cap  Z _
\rho (\hat  x)$. Since the definitions of $ \sball x \rho$ and $\locdom x
\rho$ depend on a coordinate cylinder, there may be several choices
for these sets. 
\note{ as defined above, the boundary intervals do not depend on choice of 
coordinate cyliner}
Our estimates will hold for any such choice with the
convention that when two of these objects appear in an estimate we
use the same coordinate cylinder to define both of them. 
The local domains $ \locdom x \rho$ are star-shaped Lipschitz domains
with Lipschitz constant depending only on $M$, the Lipschitz constant
for $ \Omega$. This will be helpful below as we will need to know that
various estimates hold uniformly over all local domains $\locdom x
\rho$.   See Ott and Brown \cite[p.~4376]{MR3040944} for a definition
of  star-shaped Lipschitz domains and for more details. 

We assume that the set  $D \subset \partial\Omega $ satisfies the {\em
  Ahlfors-David regularity condition}. Thus with $M$ the constant that 
controls the local behavior of our domain, 
\begin{equation}\label{AhlDav}
M^{-1}\rho \leq \sigma (\Delta _ \rho (x) ) \leq M  \rho,\  \qquad 0 < \rho  < R_0,\,
x\in  D. 
\end{equation}
Above, $ \sigma $ stands for surface measure.
We also require the sets $D$ and  $N$  to have nonempty interior
in the following quantitative sense. 
We assume that there exists $M> 0$ such that 
\begin{gather} 
 \label{DOpen} 
\text{There exists  $ x\in D$ so that 
 $\sball x  { M^ { -1} R_0}  \subset D$, }
\\
\label{NOpen} 
\text{There exists  $ x\in N$ so that 
 $\sball  x  { M^ { -1} R_0}  \subset N$. }
\end{gather}
The assumption on $N$ is used in Proposition \ref{Bogo} to solve the
equation $ \div u =f$ even 
when $f$ does not have mean value zero.   
The assumption on $D$ is used  in Appendix \ref{KornApp} to obtain
coercivity of the 
quadratic form.   Note
that this condition is only at the scale of the domain and is not
assumed to hold at every scale. In our proof of local regularity
below, we will consider mixed problems on local domains $ \locdom x \rho$ with $ \rho $ arbitrarily small. 
 We will  have freedom to specify boundary conditions on $ \locdom x \rho \setminus \partial \Omega$ and 
 will be able to guarantee that the conditions  (\ref{NOpen}) and (\ref{DOpen}) hold on all 
 local domains. 

The estimates in this paper are of two types. We will prove local
estimates for solutions that hold at scales $ \rho$  with $ 0 < \rho
< R_0$ and with a constant that depends only on $M$ and the indices
of any Lorentz spaces that appear in the estimate. In estimates over
the entire domain, the constants will also depend on the 
 collection of coordinate cylinders that
cover the boundary and such constants will be said to depend on the
Lipschitz character of $ \Omega$.  

\subsection{Function spaces} 

For $ 1 < q < \infty $, $ 1 \leq r \leq \infty$, 
we let $L^ { q,r}  $  denote the standard Lorentz space as defined in
\cite{BL:1976}, for example. 
For $ k = 1, 2, \dots$, we will use 
$ L^ { q,r}_{k}( \Omega)$ to denote   the  Lorentz-Sobolev space of functions
with $k$ derivatives in $L^ { q, r}(\Omega)$, 
\[
L^ { q,r } _k ( \Omega ) = \{ u : D^ \alpha u \in L^ { q, r} ( \Omega
) \text{ for all } |\alpha | \leq k \}
\]
with the scale-invariant norm defined by   
\[
\| u \|^q_ {L^ { q, r } _k ( \Omega ) }
=  \sum _ { |\alpha | \leq k } R_0 ^ { q(|\alpha | - k) }  \|D^ \alpha u
\| ^ q_ { L ^ { q, r } (     \Omega ) }.
\]
In the special case when $ q = r$, we will drop the second index and
note that this gives the standard Sobolev space $ L^ q _k(\Omega)$. We
will need to know that the spaces $ L^{q,r}_k ( \Omega ) $ form 
a real interpolation scale and that the spaces $
L^q_k(\Omega)$ are a 
complex interpolation scale (in $q$). See 
 \cite[Proposition 2.4]{JK:1995} for  this result in
the setting of complex interpolation and 
\cite[Proposition 2.1]{MR2763343} for real interpolation.

If $D \subset \partial   \Omega$ is a closed subset, define $ C^
\infty _D ( \bar \Omega)$ to be the collection of 
 functions that have derivatives of all orders in $ \bar \Omega$,
 whose derivatives extend continuously to $ \bar \Omega$ and 
that vanish on a neighborhood of $D$.  For $ 1 < q< \infty $, $ 1 \leq r
\leq \infty$, we define 
$L^ { q,r } _ {1,D} ( \Omega)$ to be the closure of $ C^ \infty _D(
\bar \Omega)$ in $L^ { q, r} _{1}( \Omega)$. 
Since we do not have a nice dense class in the space $ L^ { q, \infty
} _1( \Omega)$, the above definition does not serve to define the
spaces $L^ { q, \infty}_{1, D}( \Omega)$.   An alternate definition is
suggested by the work of Haller-Dintelmann, Jonsson, Knees, and
Rehberg \cite{HJKR:2012}
  who  use the 
machinery developed by Jonson and Wallin \cite{MR820626} to construct a projection $
{\cal P } : L^ q _ { 1} ( \Omega) \rightarrow L^ q _{ 1, D} (
\Omega)$.  By interpolation, we have that the map $ {\cal P} : L^ { q, r
  } _1 ( \Omega) \rightarrow L^ { q,r } _1 ( \Omega) $.  
Thus  we may  define  $ L^ { q, r } _{1, D}(
\Omega) = {\cal P } ( L^ { q, r } _ 1 ( \Omega))$ and since $ {\cal
  P}^2= {\cal P}$  the space $ L^{q,
  r}_{1,D}( \Omega)$ will be a closed subspace of  $L^ {
  q,r}_1(\Omega)$.   It is straightforward to show that the family of
spaces $L^ { q,r } _ {1, D}( \Omega)$ is an interpolation scale and
thus  with   $[\cdot, \cdot ] _{ \theta }$ denoting the operation of
complex interpolation, the 
spaces $ L^q_{1,D}( \Omega)$ satisfy
\begin{equation} 
\label{Complex}
[ L^ { q_0}_{1,D} ( \Omega), L ^ { q_ 1 }_ { 1, D} ( \Omega) ]_ \theta
= L^ { q_ \theta } _ { 1, D} ( \Omega), \quad \frac 1 { q_ \theta } =
\frac { 1- \theta } { q _0 } +  \frac \theta { q_ 1 }. 
\end{equation}
We will use  $( \cdot, \cdot ) _ {  \theta ,r } $  
to denote  the real interpolation operation and  we have 
\begin{equation}  \label{Real}
( L^ { q_0, r_0} _{1,D}( \Omega) , L^ { q_1, r_1}_{1, D} ( \Omega ) ) _ {
    \theta , r} =  L^ { q_ \theta, r}_{1,D} ( \Omega ) , \qquad
  \frac 1 { q_\theta } = \frac { 1-\theta } { q_0 } + \frac \theta { q_1} . 
\end{equation}
These results hold in our setting where the domain is Lipschitz and
the set $ D$ satisfies the Ahlfors-David regularity condition
(\ref{AhlDav}).  In fact these results hold more generally and we
refer to 
 the work of Haller-Dintelmann {\em et.~al.} \cite{HJKR:2012},
Auscher {\em et. al. } \cite{ABHR:2012} 
and 
 Brewster  {\em et.~al.}
\cite[section 6]{BMMM:2012} who consider these issues on a   larger
family of domains.

 We recall that a version of H\"older's
inequality holds in the Lorentz spaces (see \cite[Theorem 3.5]{MR0146673}). For
a set $ {\cal O }$,
\begin{equation}
\label{LorentzHolder}
\int _ {\cal O }| f g |\,d y \leq C\| f \|_ {L^ { q_0, r _ 0 } ( \cal
  O)} 
\| g\| _ { L ^ { q _ 1, r _1}( \cal O ) } , \qquad  \frac 1  { q _
  0 } + \frac 1  { q _ 1 } =1 , \ \frac 1 { r _0 } + \frac 1 { r _1} \geq
1. 
\end{equation}
Furthermore, if $1<q< \infty$ and $ 1\leq r < \infty$, then 
with $ X^*$ denoting the dual space of a Banach space $X$, we have 
\begin{equation}
\label{DualLorentz}
 L^ { q, r } ( \Omega) ^* = L^ {q', r'}( \Omega), \qquad \frac 1 q +
 \frac 1 { q' }= \frac 1 r + \frac 1 { r'} = 1. 
\end{equation}
The result (\ref{DualLorentz})  follows from a general result about duality and
interpolation, see Theorem 3.7.1 in the monograph of Bergh and
L\"ofstrom \cite{BL:1976}.

\subsection{Inequalities}

The weak formulation for the Stokes operator introduced below will use
the quadratic form 
$a : L^ { q,r } _ { 1, D} ( \Omega)
\times L^ { q', r' } _ { 1,D} ( \Omega ) \rightarrow \reals $ given by 
\begin{equation}
\label{FormDef} 
a(u,v) = 2 \int _ { \Omega } \epsilon _i ^\alpha (u) \epsilon _ i ^
  \alpha (v) \, dy .
\end{equation}
Here $ \epsilon ^ \alpha _ i (u) = \frac 1 2 ( \frac { \partial u ^
  \alpha } { \partial x _ i } +\frac  { \partial u ^ i } {\partial x _
  \alpha } ) $ denotes the symmetric part of the gradient of $u$ and
we use the convention that repeated indices are summed. 
Our conditions imply that we may find a constant $c$ depending only on
the global Lipschitz character  of $ \Omega$ such that
\begin{equation}
\label{Korn}
a(u,u) \geq c \| \nabla u \| _ { L^ 2 ( \Omega)}, \qquad u \in L^ 2 _{
  1, D} ( \Omega). 
\end{equation}
See Proposition \ref{Coercive} in Appendix \ref{KornApp}.

Next we give a small extension of a result of Bogovskii 
\cite{MR631691} 
which will
allow us to solve $\div u =f$ with $ u \in L^2_{1,D}( \Omega)$ when $f
$ is in $L^2 ( \Omega)$.
In the proof below, we will use $L^ 2 _0 ( E)$ to denote  $\{ f \in L^2 (E)
: \int _E f =0\}$.
The main fact from Bogovskii we will need is that there is 
a linear map  $B: L^2 _0( \locdom x r ) \rightarrow L^{2}_{1, \partial \locdom x r }
( \Omega)$  which satisfies $ \div B f =f$ and the norm of this
operator is bounded by $ C/r$ where the constant $C$ depends only on
the constant $M$ which controls the local  Lipschitz character of $
\Omega$. 
The argument below extends this result to general domains in a form
that is useful for the study of the mixed problem.

\note{ 
An $L^p$-result holds as well.  
 }

\begin{proposition} 
\label{Bogo} 
 Let $ \Omega $ be a connected Lipschitz domain with a
  decomposition of the boundary $ \partial \Omega = D \cup N$ and
  assume that $N$ satisfies (\ref{NOpen}).    There exists 
  a linear  map $ B : L^ 2 ( \Omega ) \rightarrow L^ 2 _ { 1, D} (
  \Omega)$ that satisfies
\[
\div B f =f 
\]
and $ \|  Bf \| _ { L^ 2 _{ 1,D}( \Omega) } \leq C \|f \|_ { L^
  2 ( \Omega) }$, where $C$ depends on the global character of $
\Omega$. 
\end{proposition}

\begin{proof}
We begin by covering $\Omega$ by a collection of local domains,
$ \locdom {x_j}{R_0/2}$, $ j =0,1,\dots,N$, and set $ \omega _j
= \locdom {x_j } {R_0}$ and $ \Omega _k = \cup _ {  j=0 } ^ k \omega
_j $, $ k=0, \dots, N$. Using the assumption (\ref{NOpen}) and
reindexing the domains $\{\omega_j \}$, we may
assume that there is a surface interval $ \sball x
s \subset \partial \omega _0 \cap  N$ with $ s $ comparable to $R_0$. 
Furthermore, since $\Omega$ is
connected, we may order the domains so that $ \locdom {x_k
}{R_0/2} \cap ( \cup _ { j = 0 } ^ { k-1} \locdom { x_j } { R_0
/2}  ) \neq \emptyset$.

 We will inductively define a sequence of maps $
P_k : L^ 2 ( \Omega) \rightarrow L^2 ( \omega _0)\times  L^2_0( \omega
_1)\times
\cdots \times  L^2 _0 ( \omega _k ) $, $k=0, \dots, N$ so that $ P_k f = ( f_0, \dots,
f_k)$  satisfies
\begin{equation}
\label{Hype}
\sum _ { j = 0 } ^ k f_ j = f \cdot \chi_{ \Omega_k } , \qquad \sum _
{ j = 0 } ^ k \| f_j \| _{L^ 2 ( \omega_j ) } \leq C \|f\|_{ L^
2( \Omega_k )}.
\end{equation}
In the first sum of (\ref{Hype}) we set $f_j$ to be zero outside
$ \omega_j$. We define $ P_0$ by $P_0 f =f \cdot  \chi _ { \omega_0 } $ and
clearly (\ref{Hype}) holds for $k=0$. 

Before we define $ P_k$, $k =1, \dots,N$, we observe that since
$ \Omega _ { R_0 /2 } ( x _k ) \cap ( \cup _ { j = 0 } ^ { k-1} \Omega
_{R_0/2} ( x_j )) \neq \emptyset$, we have that the measure of $\omega
_k \cap \Omega_{k-1}$ is comparable to $R_0^2$. If $ f \in
L^2 ( \omega _k )$, we write $ f\cdot \chi _ { \Omega _k }  = g  + f _k$
with 
\[
f_k = f \cdot \chi _{ \omega _k \setminus \Omega _{ k-1}} - \frac 
{\chi_{ \omega _k  \cap \Omega _ { k
-1 }}}
 {
|\omega _k \cap \Omega _{ k -1} | } 
\int _ { \omega _k \setminus \Omega _{k-1} } f \, dy 
\]
and $ g =f - f_k$. The inequality of Cauchy-Schwarz and our
observation that the measure of  $ \omega _k \cap \Omega_{k-1}$  is
comparable to $R_0^2$ implies that 
\[
\| f_k \|_{ L^2 ( \omega_k )} + \|g \|_{ L^2 ( \Omega _{ k -1})}
\leq \| f\| _ { L^2 ( \Omega_k)}.
\]
We use our induction hypothesis that (\ref{Hype}) holds for $ k
-1$ and let $(f_0, \dots, f_{ k-1}) = P_{k-1} g$. It is easy to
see that $ P_k f = ( f_0, \dots, f_k )$ is the desired decomposition. 

Given the decomposition of $f \in L^2 ( \Omega)$, $f=(f_0, \dots,
f_N)$, we may use the result of Bogovskii  in each $ \omega _j$, $j
=1, \dots N$,
to find $u_j$ which satisfies $ \div u _j =
f _j$ and $ u _ j \in L^ 2 _ { 1, \partial \omega _j }( \omega _j
)$. Recall
that $L^ 2 _ { 1, \partial \omega _j }( \omega _j ) $ is 
the Sobolev space of functions with one derivative in $L^2$ and which
vanish on the boundary $ \partial \omega_j$. To define $u _0$, we use the
existence of the surface interval $\sball x s$ noted above to find a
vector-valued
function $ \eta \in C^\infty _D(\bar \Omega)$,
$\supp \eta \subset \omega _0 \cup \sball x s $ and so that 
$ 1 = \int _ {\sball x s} \eta \cdot \nu\, d\sigma = \int _ { \omega
_0 } \div \eta \, dy. $  We let $ u _ 0 = v _0 + \eta \int _ {\omega
_0 } f_0 \, dy $ where $ v _0 $ is the solution of $ \div v _0 = f_0
- \div \eta \int _{ \omega _0 } f_0\, dy $ guaranteed by Bogovskii's
result. Finally, we let $ u = \sum _ { j =0 } ^ N u _j$ 
where we have defined each $u _j $ to be zero outside
$ \omega_j$. 
It is immediate
that $ u \in L^2 _ { 1, D} ( \Omega )$ and satisfies $ \div u
=f$. 
\end{proof}

\section{Weak solutions}
\label{WeakSolutions}
The goal of this section is to show that we can solve the mixed
problem when the right-hand side of \eqref{WMP} lies $L^ { 2,1} _{ 1,D}
( \Omega ) ^*$, the dual of the Lorentz-Sobolev space $ L^ { 2,1 } _ { 1,D}
(\Omega) $.   Since functions in  $L^ { 2,1 } _ 1 ( \Omega)$ are
continuous, the Dirac delta measure lies in the dual of $L^ { 2, 1} _1
( \Omega) $ and thus we are able to construct the Green function as
the solution of (\ref{MP}) when $f$ is the Dirac delta measure,
$f_N=0$ and $g =0$. 

We give a weak formulation of the mixed problem. We will
need to consider this formulation, not only on the Sobolev space $L^ 2
_ { 1, D} ( \Omega) $ but also on Lorentz-Sobolev spaces with indices
$q$ near 2. 
To simplify the notation below, we introduce spaces $ S_ {q,r} = L ^ {
  q, r } _ { 1,D} ( \Omega) \times L^ { q,r } ( \Omega)$ and  then  $ S^
* _{ q, r} $ will  denote the dual of $S_ { q, r}$. As with the Lorentz and
Lebesgue spaces, we will drop the second index when both indices take
the same value and use $ S_q = S_ {q,q}$. 

\note{ The particular form is not important, provided we have coercivity
  and boundedness. 
}

We give a weak formulation of (\ref{MP}) in the space $S_{
  q,r}$. For this problem, we will require 
$ f\in L^ { q, r} _{ 1, D}(\Omega) ^*$, 
$g \in L^ { q,r}( \Omega) $ and $f_N$ to lie in the dual
of the image of $L^ { q, r} _{1, D} ( \Omega)$ under the trace
operator.  We give a weak formulation of (\ref{MP})  in the special
case that $f_D=0$. 
We say that $(u, p) $ is a weak solution of (\ref{MP}) with $f_D=0$ if we have
\begin{equation}
\label{WMP}
\left\{
\begin{array}{ll}
 a(u,\phi) - \int _ { \Omega } p \div \phi  \, dy 
\\ \qquad \qquad = \langle f, \phi\rangle
+ \langle f_N, \phi \rangle _{ \partial \Omega}
-
 \int _ \Omega \nabla g \cdot  \phi \, dy,
\qquad 
& \phi 
\in L^ { q' , r'}_{1,D} ( \Omega)  \\
 - \div u = g, \\
( u,p) \in S_ {q,r}.
\end{array}
\right.
\end{equation}
\note{ 
It is not clear if there are simple conditions on $g$  that guarantee the
right-hand side is in the dual of $L^{ q',r'}_{1,D} ( \Omega)$. 
This weak formulation is not quite right--we weren't allowed to
  integrate by parts to obtain $ \int g \div \phi$. When $g$ is zero,
  this is not a problem, but for general $g$, there is something more
  to be said. 

This might affect the Caccioppoli inequality, too.
}
The three terms on the right of the first equation in (\ref{WMP}) have different roles in the 
 Stokes system. However,  for the purposes of
establishing existence of solutions, it 
makes sense to simplify the problem by treating 
the three terms as  one.  We will consider the case when the
right-hand side of the first equation in (\ref{WMP}) gives an element
in  the dual of $L^ { q', r'} _{ 1,
  D} ( \Omega)$.   It is not obvious how to give general conditions on
$ f_N$, $f$, and $g$ which guarantee that this will happen. However,
our construction of the Green function will not require us to consider
this question.

\note{  Can we find an example of a function $g$ in $L^2$ for which the right-hand side is not in $L^2_{-1}$? }

We introduce a map $ T: S_ {q,r} \rightarrow S_ {q',r'}^*$ given
by $ T(u,p ) = ( \lambda , \mu)$ where  $ \lambda $ is in the dual of
$ L^ { q', r' } _ { 1, D} ( \Omega)$ and $ \mu$ is in the dual of $ L^ {
  q', r'} ( \Omega)$ and are given by 
\begin{align} 
\label{Tdef1}
\lambda ( \phi ) &= a(u,\phi) -\int_\Omega p \div \phi \, dy  \\
\label{Tdef2}
\mu(h)  & = - \int _ \Omega h \div u \, dy .
\end{align}
 Thus, establishing the existence of a solution to
(\ref{WMP}) is equivalent to showing that the map 
\begin{equation}
\label{UpsideDown}
T : S_ {q, r}  \rightarrow S_ {q',r' } ^* \text{ is invertible.}
\end{equation}
The next theorem outlines the main steps in establishing
(\ref{UpsideDown}). 
\begin{theorem} 
\label{Exist}
Let $ \Omega $ be a Lipschitz domain and suppose that the
decomposition 
of the boundary $ \partial \Omega = D\cup N$ satisfies
(\ref{AhlDav}), (\ref{DOpen}), and  
(\ref{NOpen}). Then $ T$ is a bounded operator satisfying 
\begin{align}
& T : S_ {q,r } \rightarrow S_ {q',r'}^*, \qquad 1 < q <  \infty , 1 \leq
  r \leq \infty \label{Bounded} \\
& T : S_ 2\rightarrow S_2 ^* \text{ is invertible}. \label{Energy}
\end{align}
Moreover, there exists $\kappa > 0$ such that 
\begin{align}
& T: S_q \rightarrow S_{q'} ^ * \text{ is invertible if } |\frac 1 2 -
  \frac 1 q | < \kappa \label{Jiggle} \\
& T ^ { -1 } : S _ { q' ,r' } ^ * \rightarrow S_ { q, r } , \ |\frac
  1 2 - \frac 1 q | < \kappa, \ 1 < r \leq \infty. \label{Useful}
\end{align}
The constant $ \kappa$ depends only on the Lipschitz character of $
\Omega$. The norms of the operators depend on the Lipschitz character
of $ \Omega$ and the exponents $q$ and $r$ for the Lorentz spaces. 
\end{theorem}

\begin{proof} The first claim (\ref{Bounded}) is an easy consequnce of
the extension of H\"older's inequality to the Lorentz spaces, 
(\ref{LorentzHolder}). The second claim (\ref{Energy})
 is a standard  existence
theorem for weak solutions of the Stokes system, but  we sketch the details due to 
its fundamental importance. We follow the argument given by
Maz'ya and Rossman for polyhedral domains \cite{MR2563641}. 

In the proof of the second claim (\ref{Energy}), we
will work 
 in  the subspace of functions that
are divergence free. 
Thus we let
\[
H = \{ u \in L^ { 2} _{ 1, D} ( \Omega) : \div u =0\}
\]
and  we define $ H^ \perp$ to be the orthogonal complement of $ H$
in $L^ 2 _ { 1, D}( \Omega) $.
\comment{ where we use
\[
\int _ \Omega \frac { \partial u ^ \alpha} {\partial x _i } \frac {
  \partial v ^ \alpha }{\partial x _ i }\, dx 
\]
as the inner product 
on $L^2 _ {1, D} (
\Omega)$. } 
  According to Proposition \ref{Bogo}, we have a map $ B' :
L^ 2 ( \Omega ) \rightarrow L^ 2_ { 1, D} ( \Omega)$ with $\div B'f =
f$. We let $ P $ be the orthogonal projection from $L^ 2 _{1, D}(
\Omega)$ 
onto $H ^ \perp$ and
then define $ B = PB'$. It is easy to see that $ B : L^ 2 ( \Omega )
\rightarrow H ^ \perp$ is an isomorphism and $ \div Bf =f$. 

Given $ (\lambda, \mu )$ in $S_2^*$, we need to find
$(u,p) $ in $ S_2 $ which satisfies $T(u,p) =
(\lambda , \mu )$. 
Since  $ \mu $ lies in $L^2 ( \Omega)^*$, 
 there exists $ g\in L^2 ( \Omega)$ such that 
$ \mu ( h) = \int_ { \Omega } gh\,d y $.
From  Proposition \ref{Bogo}, we may find 
$ w = Bg$ in $L^ 2 _{1, D} (\Omega)$ which satisfies $ \div w =g$.  
The form $a$ is coercive  on $ L^ { 2} _{1, D}(
\Omega)$ (see (\ref{Korn})), 
hence  it is also coercive on the subspace $H$ and by
Lax-Milgram we may find 
$ v \in L^ { 2 } _{1, D}( \Omega)$ which satisfies 
\begin{equation}
\label{ExA}
a(v,\phi) = \lambda ( \phi) - a( w, \phi), \qquad  \phi \in H.
\end{equation}
We let $u = v + w$ and turn to constructing the pressure $p$. As a step in this
direction, we let $ \tau (h ) = \lambda ( Bh) - a ( u, Bh) $ for $ h
\in L^ 2( \Omega)$. As $ B : L^ 2 ( \Omega ) \rightarrow L^ 2 _ { 1,
  D} ( \Omega)$ is continuous, it follows that $ \tau \in L^ 2 (
\Omega) ^*$ and hence we may find $ p \in L^ 2 ( \Omega)$ so that 
\[
\tau(h) = - \int_\Omega ph\,d y.
\]
Since $B ( L^2 ( \Omega)) = H ^ \perp$, we have
\begin{equation}
\label{ExB}
a(u, \phi) - \int_\Omega p \div \phi \, dy = \lambda ( \phi), \qquad \phi \in
H^ \perp.
\end{equation}
Combining (\ref{ExA}) and (\ref{ExB}) and observing that $ \div u =
g$, we have found a pair $(u,p)$ which satisfies $ T(u,p) = (\lambda ,
\mu)$. 

To establish that $ T $ is injective, we suppose that $T(u,p)
=(0,0)$ with $(u,p) \in S_2$. If $ \mu$ as defined in (\ref{Tdef2}) is
zero, we have $ \div u =0$. Recalling  (\ref{Tdef1})  it follows that $
a(u,u) =0$ and then the 
coercivity of $a$ implies that $ u =0$. Once we have that $u=0$, 
we  let $ \phi = B p$ in the first line of (\ref{WMP}) and conclude
that $ \int_ \Omega p^2  \, dy =0$ and hence $p =0$. 

We turn to the third statement in Theorem \ref{Exist}, (\ref{Jiggle}).  The
family $S_q$ is a complex interpolation scale and by Corollary 4.5.2 in
Bergh and L\"ofstrom \cite{BL:1976} it follows that $S_q^*$ is also a
complex interpolation scale. Using (\ref{Bounded}) and (\ref{Energy}),
(\ref{Jiggle}) follows from a general result of Sneiberg
\cite{MR0634681} (see also Tabacco-Vignati and Vignati
\cite{TV:1988}).

If we choose $ q$ with $ \frac 1 2 < \frac 1 q < \frac 1 2 + \kappa
$ with $ \kappa $ as (\ref{Jiggle}), we obtain 
\[
T^ { -1} : ( S^*_q , S^ * _ {q'}) _ { 1/2, r} \rightarrow  (S_{q'} , S_
q )_{ 1/2, r}.
\]
At least for $ 1<r \leq \infty$, we have $( S^*_q , S^ * _ {q'}) _ {
\theta, r} = ((S_ {q'}, S_q)_{\theta ,r'})^*$ (see \cite[Theorem
3.7.1]{BL:1976}).  Thus using our characterization of the real
interpolation spaces for $L^q_{1, D}( \Omega)$ in (\ref{Real}) we obtain
\eqref{Useful}.  

\note{ 
Is the projection needed in the definition of $B$? 

Check that $ B: L^2 (\Omega) \rightarrow H^ \perp$ really is onto. 

} 
\end{proof}

\section{Local regularity} 
\label{Local}
In this section, we will need the following version of the Poincar\'e
inequality. A proof may be found in  our previous work
\cite[Appendix A]{MR3040944}. 
  Below,  we let $\average _ { \locdom x r } u \, dy = 
|\locdom x  r|^{-1} \int_{\locdom x r} u\, dy$ denote the average of
the function $u$ on $ \locdom x r $. 

Let $ u $ be in $L^ q _ {1, D} ( \Omega)$ and let $ \bar u _ {x, r}$
be defined by 
\[
\bar u _ {x, r } = \left \{ \begin{aligned}  0 , \qquad & \dist (
  \locdom x r , D ) = 0\\
\average _ { \locdom x r } u \, dy , \qquad & \dist ( \locdom x {  r }
, D) > 0 . 
\end{aligned}
\right. 
\]
Then for $ 1 \leq q < 2 $ we have 
\begin{equation}
\label{PoSo}
\left (\int  _ { \locdom x { r } } |u- \bar u _ { x, r}| ^ { 2q / ( 2-q)} \, dy \right )
^ { \frac  1 q - \frac 1 2 } \leq C \left ( \int _ { \locdom x { 2r} }
|\nabla u | ^ q \, dy \right ) ^ { 1/q}.
\end{equation}
This assumes that the set $D$ satisfies condition (\ref{AhlDav}). 
Note that the  integral on the right-hand side of (\ref{PoSo}) is
over a larger set.  In the case where $ \dist(\locdom x r , D) =0$,
the Ahlfors-David condition (\ref{AhlDav})  guarantees that $u$
vanishes on a subset of
$\partial \locdom x {2r}$ with measure proportional to $r$ and this
allows us to obtain an estimate for $u$ on $\locdom x {2r}$. 
 The expansion of the  domain of integration is
not needed  in the case  where $ \bar u _ {x, r } \neq 0$.

The existence of the Green function is an immediate consequence of
Theorem \ref{Exist} since we may define $ G(x,\cdot) ,\Pi (x, \cdot) $
by
$(G^ { \alpha \cdot } (x, \cdot), \Pi^ \alpha (x, \cdot) ) =
T^ { -1} ( e_ \alpha\delta _x , 0 )$. However, more work is needed to obtain
the estimates of our main result, Theorem \ref{Green}. Thus we
consider the local regularity of solutions of the boundary value
problem (\ref{WMP}) before giving the argument for existence.

 \begin{proposition} 
\label{Reverse}
Let $ \Omega$ be a Lipschitz domain with $\partial\Omega = D \cup N$ and
$D$ satisfying
(\ref{AhlDav}).  Suppose that $(u, p ) \in S_{q_1}$ is a solution of
the weak mixed problem with $ f_N =0$ and $g =0$.  Let $ \locdom x
\rho$ be a local domain and suppose that $ \eta $ is a smooth cutoff function
which is one on $\locdom x \rho$ and zero outside $\locdom x {
2\rho}$. There exists a positive number $ \kappa$ so that if $
\frac 1 { q_1} < \frac 1 2 + \kappa$, then for each $q$ with $
\frac 1 2 > \frac 1 q > \frac 1 2 - \kappa$, we have
\begin{multline*}
\| p \| _ {L^ q ( \locdom x \rho ) }+
\| \nabla u \| _ {L^q ( \locdom x \rho )} 
\\
\leq \| \eta f \| _ { L^ {q'} _{1,D} ( \locdom x {2 \rho} ) ^*} +
\frac 1 \rho  (  \| \nabla u \| _ { L^ {\tilde q }  (
  \locdom x { 4 \rho} ) }  +\|p\|_ {L^ {\tilde q} ( \locdom x { 2\rho})}) .
\end{multline*}
Here, we define $ \tilde q $  by $ \frac 1 {\tilde q} =  \frac 1 2 + \frac 1
  q $. 
 \end{proposition}
\note{
 It might make more sense to apply the norm in $L^ {q'} _{ 1, D'} (
\locdom x { 2 \rho} )^* $ to $\eta f$, however in most of our
applications,  $ \eta f $ is zero
and thus there is no point in going to the extra trouble. 
}

\begin{proof} We fix a local domain $ \locdom x \rho$ and cutoff
  function $ \eta $ as in the statement of the Proposition. We will
  show that $(\eta (u - \bar u _ {x, 2 \rho}), \eta p)$ is a solution of a mixed problem in $
  \locdom x { 2\rho}$.  We note that our definitions of $ \bar u _ {x,
    2\rho}$ and of $ \eta$ guarantee that $ \eta ( u - \bar u _ { x,
    2\rho})$
will vanish on $D$.  To define the Dirichlet set $ D' \subset
\partial  \locdom x { 2 \rho}$, we will consider two cases. Case 1: $\bar
  \Omega _ { 2 \rho }(x)  \subset \Omega$. In this case,  $\locdom x{2\rho}$
  is a disk and  we define $ D'$
  to be an arc of length $ \pi \rho$ in  the boundary. Case 2: $\partial
  \locdom x { 2\rho } \cap \partial \Omega \neq \emptyset$. In this
  case, we let  $ {\cal S}  = \partial \locdom x {2 \rho}\cap \Omega \cap \{ y : |x_1 -
  y _1 | = 2\rho\}$ denote the sides of $ \locdom x { 2\rho}$ and then
  we put $ D' = (D \cap  \partial \locdom x {2\rho}) \cup { \cal S }$.
Finally, we set $ N'  =
  \partial  \locdom x { 2\rho}\setminus D'$ and 
 leave it as an
  exercise to check that the decomposition $ \partial \locdom x
  {2\rho} = D'\cup N' $ satisfies the Ahlfors-David
  regularity condition (\ref{AhlDav}) and the conditions  (\ref{DOpen}), and
  (\ref{NOpen}) with constants that are independent of 
$x$ and $ \rho$. 
Thus, we may find $ \kappa$ that is independent of $x$ and $\rho$ 
 so  that (\ref{Jiggle}) of  Theorem \ref{Exist} holds 
for $q$ with $ | \frac 1 2 - \frac 1 q | < \kappa$. 

We will show that, with $T$ the map defined in 
(\ref{Tdef1}-\ref{Tdef2}) for $ \locdom x {2\rho}$, we have  
\begin{multline}
\label{ReverseEstimate}
\| T(\eta (u - \bar u _ { x, 2\rho}), \eta p )\| _ { S_{q'} ^ *
  (\locdom x { 2\rho} ) } \\ 
\leq \frac C \rho ( \| \nabla u  \| _ { L^ { \tilde   q }  ( \locdom x
   {4 \rho}) } +  
\|p  \| _ { L^ { \tilde   q }  ( \locdom x {2 \rho}) } )
+ \| \eta f \| _ { L^ { q' } _ { 1, D' } ( \locdom x { 2\rho} ) ^*}. 
\end{multline}
With this, the estimate of the theorem will follow from the definition
of $ \eta$ and Theorem \ref{Exist}. Our assumption that $ (u,p) \in S_
{q_ 1} (\Omega)$ with $ \frac 1 { q _ 1} < \frac 1 2 + \kappa$
guarantees that $ (\eta( u- \bar u _{x, 2\rho}) , \eta p) $ is the unique solution (\ref{WMP})
in $S_{q_1}( \locdom x { 2 \rho})$. 

Thus, we turn to the proof of (\ref{ReverseEstimate}).  We fix $ q$
with 
$\frac 1 2 > \frac 1 q > \frac 1  2  - \kappa$ and choose 
$\phi \in L^ { q' } _ { 1, D' } ( \locdom x { 2 \rho } ) $.  We use $
\eta \phi$ in the weak formulation of the mixed problem in $\Omega$ satisfied by $
u$  to obtain 
\[
a( u- \bar u _ { x, 2\rho} , \eta \phi ) + \int_ \locdom x { 2\rho} p \div( \eta \phi )\, dy =
\langle f,  \eta \phi \rangle .
\]
Using the product rule, we may rewrite this as 
\begin{multline} 
\label{Happy}
a ( \eta ( u - \bar u _ { x, 2\rho} ), \phi ) + \int_ { \locdom x { 2 \rho}}
\eta p \div \phi\,d y 
= \langle \eta f , \phi \rangle    \\
+ \int_ { \locdom x { 2\rho} } - p \nabla \eta \cdot  \phi  
 + \epsilon^ \alpha _ i ( u - \bar u
_ {x , 2\rho } ) \frac {\partial \eta }{ \partial y _ i  } \phi ^ \alpha
- \frac { \partial \eta }{ \partial y _ i } ( u - \bar u _ { x , 2\rho })
^ \alpha \epsilon _ i ^ \alpha ( \phi ) \, dy \\
= \langle \eta f , \phi \rangle  + I + II + III.   
\end{multline}
We claim that 
\begin{equation}
\label{Grumpy} 
|I+II+III| \leq \frac 
 C  \rho (  \| p \| _ { L^ {\tilde q}  ( \locdom x { 2 \rho}
  ) }
+ \| \nabla u \|_{ L^ {\tilde q }  ( \locdom x {
    2 \rho} )})  
\|\nabla \phi \| _      { L^ {q'} ( \locdom x { 2\rho} ) }.
\end{equation}
Furthermore,  it is clear from the product rule and the
Sobolev-Poincar\'e inequality (\ref{PoSo}) that $ - \div ( \eta ( u
- \avg u x { 2\rho}) )$  satisfies 
\begin{equation}
\label{Sneezy}
\| 
 - \div ( \eta ( u
- \avg u x { 2\rho} ) )\| _ { L^q ( \locdom x { 2 \rho }) } 
 \leq 
\frac C \rho \| \nabla u \| _ { L^ { \tilde q } ( \locdom x { 4 \rho}
  )}.
\end{equation}
Combining (\ref{Grumpy}) and  (\ref{Sneezy}) gives
(\ref{ReverseEstimate}) and hence the Proposition. Note that we may
replace $D'$ in (\ref{ReverseEstimate}) by $D$ in the  Proposition
since the norm of $ \eta f$ in  $L^{p,1}_D(\Omega)^*$ is an increasing function
in the set $D$. 

Thus it remains to prove the claim (\ref{Grumpy}). 
To estimate  $I$ and $II$,  we use  H\"older's inequality and the
Sobolev-Poincar\'e inequality (\ref{PoSo})
\begin{align*}
| I | +|II|
 & \leq \frac C  \rho (  \| p \| _ { L^ {\tilde q}  ( \locdom x { 2 \rho}
  ) }
+ \| \nabla u \|_{ L^ {\tilde q }  ( \locdom x {
     2 \rho} )}) 
 \| \phi \| _      { L^ {\tilde q '} ( \locdom x { 2 \rho} ) }
\\
&\leq 
\frac C  \rho (  \| p \| _ { L^ {\tilde q}  ( \locdom x { 2 \rho}
  ) }
+ \| \nabla u \|_{ L^ {\tilde q }  ( \locdom x {
    2 \rho} )}) 
\|\nabla \phi \| _      { L^ {q'} ( \locdom x { 2\rho} ) }.
\end{align*}
To estimate $III$, we 
use H\"older's inequality and then apply  the Sobolev-Poincar\'e  
inequality \eqref{PoSo}
to $ u - \bar u _ {x , 2 \rho }$ to obtain 
\begin{align*}
|III| &\leq \frac  C \rho
\| u - \bar u _ { x, 2\rho }\| _ {L ^ {q}( \locdom x { 2 \rho} ) } 
\| \nabla \phi \| _ {L^ {q'} ( \locdom x { 2\rho } )}
\\
&\leq \frac C \rho   \| \nabla u \| _ { L^ { \tilde q } ( \locdom x { 4\rho} ) }
 \| \nabla \phi \| _ { L^ {q'} ( \locdom x { 2\rho} )} .\\
\end{align*}
Combining the last two displayed  estimates gives the claim
(\ref{Grumpy}) and hence the Proposition. 
\end{proof}

We will need the following version of the Caccioppoli inequality  from 
 Giaquinta (\cite[Theorem 1.1]{MR641818}).  Let $ u$ be a
solution of (\ref{WMP}) with $ f = 0$ and $ g =0 $  on $\locdom x { 2
  \rho}$
and $ f_N = 0 $ on $\sball  x { 2\rho} \cap N$, then we have 
\begin{equation}
\label{Caccioppoli} 
(\average _ { \locdom x \rho } |\nabla u |^2 \, dy ) ^ { 1/2}  \leq 
C \average _ { \locdom x { 2\rho }} |u|\, dy .
\end{equation}
Giaquinta and Modica prove this estimate with an $L^2$ norm on the
right and in the interior. The argument may be adapted to work near
the boundary. It is well-known that we may vary the exponent on the
right-hand side, see \cite{MR1239172} for example.

\begin{corollary} 
\label{Continuous}
Suppose that $ \Omega$ is Lipschitz domain  and $ \partial \Omega = D
\cup N$ with $D$ satisfying 
(\ref{AhlDav}).  Let $ \kappa $ be as in Proposition \ref{Reverse}. 
Suppose that $ (u,p) \in S_{q}$, $ 1/q < 1/2 + \kappa$  is a solution of (\ref{WMP}) in $\Omega$ with $ f_N
=0$,    $ f_ D=0$, and $g=0$  and $ f $ is zero on $\locdom x {2\rho}$. Then we have
the estimates 
\begin{gather}
\label{LocalHolder}
|u(z) - u (y) | \leq C \left ( \frac { |z-y|} \rho \right) ^ \gamma 
 \average _ { \locdom x {2\rho}} |u(y)- \avg u x {2\rho} |\, dy , \qquad y, z \in \locdom x
\rho 
\\
\label{MVT}
|u(x) | \leq C  \average _ { \locdom x \rho} |u(y) |\, dy .
\end{gather}
\end{corollary}

\note{ Need to check hypotheses needed for caccioppoli. } 

\note{ A comment about exponent $\gamma$. }

\begin{proof} The first estimate  (\ref{LocalHolder}) follows from
Proposition \ref{Reverse},    Morrey's inequality, and the Caccioppoli
inequality (\ref{Caccioppoli}). This gives the H\"older estimate with
$0< \gamma < 2 \kappa$ with $ \kappa$ as in  Proposition \ref{Reverse}.   
To
obtain the second inequality, we write
$
|u(x) | \leq |u(x) - u (z) | + |u(z) |. 
$
We may average this inequality with respect to $z$ in $ \locdom x {\rho/2}$
and then use   (\ref{LocalHolder}).   
\end{proof}

\section{Proof of Theorem \ref{Green}}
\label{TheEnd}

We are ready to give the proof of our main result, 
Theorem \ref{Green}, from the introduction. 
\begin{proof}[Proof of Theorem \ref{Green}]
We begin by noting that functions in the Lorentz-Sobolev space
$L^ { 2,1 } _{1,D} ( \Omega)$   are continuous. To see this, we
observe that since $ \Omega$ is Lipschitz,  we have an extension
operator 
 $E: L^{ 2,1} _1( \Omega) \rightarrow 
L^ { 2,1 } _1( \reals ^2)$. If $f$ is 
 in $L^{ 2,1}_1( \reals^2)$, we may represent $f$ as a first order
 Riesz potential 
and use H\"older's inequality  in Lorentz spaces (\ref{LorentzHolder})
to conclude  that  $f$ is bounded and continuous
(see  \cite{MR607898}, for example).  
Thus the  Dirac delta measure concentrated at a
point $x$ in $ \Omega $ lies in $ L^ { 2, 1 } _ { 1, D} ( \Omega)
^*$.  For each $ x \in \Omega$, $ \alpha =1$ or 2, we put 
\[
( G^{ \alpha \cdot} ( x, \cdot  ) , \Pi^ \alpha (x, \cdot) ) 
= T^ { -1 } (e_\alpha \delta _x , 0), 
\]
where  $e_\alpha$ is a unit vector and $ T $ is the map defined
in (\ref{Tdef1}-\ref{Tdef2}).  From  (\ref{Useful}), we have that  $
T:S_{2,\infty} \rightarrow S_{2,1}^*$ is invertible and hence
\[
\| G(x, \cdot ) \| _ {L^ { 2, \infty }_{1,D} (\Omega)} + \| \Pi(x,
\cdot)  \| _ { L ^ { 2, \infty } ( \Omega) }  \leq C.
\]
Thus, we have established (\ref{GreenLorentz}) of Theorem
\ref{Green}. 
Since we have  $ \nabla G(x,\cdot)$ in
$L^{2, \infty} (\Omega) $ it follows that  for each local domain $
\locdom x \rho$
\[
\int _ { \locdom x \rho } |\nabla _y G(x,y) | \, dy \leq C\rho. 
\]
Using the Sobolev-Poincar\'e inequality (\ref{PoSo}), we obtain
that $ G (x, \cdot )$ is in the space $BMO_D(\Omega)$ as defined by
Taylor, Kim, and Brown \cite[p.~1579]{TKB:2012}. Following the arguments 
from Taylor, Kim, and Brown,   the local 
regularity estimates in Corollary \ref{Continuous} imply  the upper bound
(\ref{GreenUp}) and the H\"older estimate (\ref{GreenHolder}) for $
G(x, \cdot)$.  
The estimate for the gradient of $G $ in (\ref{GreenGrad}) follows
immediately from Proposition \ref{Reverse}.

To establish the representation formula, (\ref{GreenRep}), we observe
that if $u$ is a solution of (\ref{WMP}) with $f_N \in L^t(N)$, $f \in
L^ t ( \Omega)$,  and  $ \nabla g \in L^ { t } ( \Omega)$ for some $t>1$,
then we have that the solution $(u,p)$ lies in $S_q$ for some $q>2$.  From
the weak formulation (\ref{WMP}) for the equation for $G$, we obtain 
\[
a( G^ { \alpha \cdot} (x, \cdot) , u ) - \int _ \Omega \Pi ^ \alpha (
x, y) g(y) \,dy   = u^ \alpha (x) - \int _ \Omega \nabla g \cdot G^ {
  \alpha \cdot }(x,y) \, dy  .
\]
On the other hand, using $ G^ { \alpha \cdot} (x, \cdot) $ as the test function in the
weak formulation of the boundary value problem satisfied by $u$ gives 
\[
a(u, G^ { \alpha \cdot } (x, \cdot ) ) = \int _ \Omega G ^ {\alpha
\beta} f^ \beta (y) \, dy + \int _ { N } G^ { \alpha \beta } (x, y)
f_N^ \beta (y)  \, d\sigma .
\]
We recall that the form $a(\cdot, \cdot)$ is symmetric and 
equate  the two expressions for $a(u, G^ { \alpha \cdot} (x, \cdot)
)$ which  gives the representation formula (\ref{GreenRep}).   If $ g $ is
smooth and compactly supported in $ \Omega$, then integration by parts
shows that the 
term in (\ref{GreenRep}) involving
$G$ and derivatives of $g$ is zero. This
allows one to see that $G$ is a Green function as defined 
before Theorem \ref{Green}.

Next we prove the symmetry property for the Green function. On a formal
level, the proof follows quickly from the weak formulations of the
problems for $ G ^{ \alpha \cdot } (x, \cdot)$ and $ G^{ \beta
  \cdot}(y, \cdot) $ and the symmetry of the form $a(\cdot, \cdot)$. 
These combined give
\[
a( G^ { \alpha \cdot } (x, \cdot) , G^ { \beta \cdot}(y,\cdot)  )
= G^ { \alpha\beta} (y,x) = G^{ \beta \alpha } (y,x).
\]
However, this formal argument fails since we do not have that $ G(x,
\cdot)$ lies in $L^{ 2,1 } _{1, D}( \Omega)$.  The argument can be
made rigorous by approximating the $ \delta $-functions  by
normalized characteristic functions of balls centered at $x$ and $y$
and letting the radii tend to zero. See the work of Taylor, Kim, and
Brown \cite{TKB:2012} for a similar argument. 
\end{proof}

We close with a few related questions for further investigation. Many
of the 
methods used in the work reported here apply also to variable
coefficient equations. 
It would be interesting
to extend the construction of the Green function to operators with
non-constant coefficients that are 
 modelled on  the Stokes system and to domains that are more
general than Lipschitz. 

One point that remains open is finding a representation for the
pressure. One may define another row of the Green matrix to be a
solution of the system 
\[
- \Delta G^{3\cdot}(x, \cdot ) + \nabla \Pi^3(x, \cdot)  = 0 , \qquad 
-\div G^ { 3\cdot}(x, \cdot) = \delta _x
\]
and use $G^ { 3 \cdot}$ to give a representation of the
pressure. See, for example, the work of  Maz'ya and Rossmann
\cite{MR:2005}.  One expects that $ G^ {3\cdot} (x, \cdot ) $ lies in
$L^{2, 
  \infty}( \Omega)$ uniformly for $ x\in \Omega$. However, we are
unable to obtain such estimates.

\note{
Need to prove symmetry. This may require a limiting argument. 

Representation for the pressure. Do we have the right bvp?

} 

\appendix
\section{Korn inequality}
\renewcommand*{\thesection}{\Alph{section}}
\renewcommand{\theequation}{\Alph{section}.\arabic{equation}}
\label{KornApp}

In this appendix, we 
provide a proof of  the Korn inequality 
\eqref{Korn}.   While the inequality is well-known, 
the standard proof (see \cite[Theorem 2.7]{MR1195131}, for example) uses a
contradiction argument and does not seem to give any information about
the relationship  between the constant and the geometry of the domain. 
Since we did not find a proof that makes clear the relation, we
include a short argument exhibiting the dependence of the constant on the geometry of the domain. 
The argument uses duality and the 
 result of Bogovskii \cite{MR631691} as given in    Proposition
 \ref{Bogo}.    This line of argument 
 is  well-known and  an application to proving a different version of the Korn
 inequality can be found  in 
\cite{RD:2012}. 

\note{Definition of star-shaped convex domain is not given.  
Need $D$ to contain an open interval  to obtain coercivity.  } 
\begin{proposition}
\label{Coercive} 
 Suppose that $ \Omega  $ is a Lipschitz domain 
with a decomposition of the boundary $ \partial \Omega = D \cup N$ and
that $ D$   satisfies (\ref{DOpen}). 
Under these assumptions, we have the estimate 
\[
\| \nabla u \| _ { L^ 2 ( \Omega) } \leq C \| \epsilon (u) \| _{L^2 (
  \Omega) }, \qquad u \in L^ 2 _ { 1, D}( \Omega )
\]
with a constant $C$ that depends only the dimension $n$ and $M$. 
\end{proposition}

\begin{remark}  In our application it is important to note that the
  constant in this estimate may be taken to hold uniformly over all
  local domains $ \locdom x \rho$ provided the boundary conditions
  satisfy (\ref{DOpen}). 
\end{remark}

\begin{proof} The  estimate of this Proposition is dual to the result of Proposition
  \ref{Bogo}.  We begin by observing that it suffices to prove the
  estimate for $u \in C^ \infty _ D(\bar \Omega )$. We let $ \mu^\alpha _ {
    i} (u) = \frac 1 2 ( \frac { \partial u ^ \alpha} {\partial
    x _ i } - \frac { \partial u ^ i } { \partial x _ \alpha })$
  denote the anti-symmetric part of the gradient. The proposition will
  follow once we show 
\[
\int _ \Omega  | \mu  (u) |^2 \, dy
\leq C \int _ \Omega | \epsilon(u) |^2 \, dy. 
\]
Thus, we fix $i$ and $ \alpha$,  choose $f$ in $L^ 2 ( \Omega)$, 
and use 
Proposition \ref{Bogo} to  find $v \in L^ 2 _ { 1} ( \Omega)$ which is  a
solution of  $ \div v = f$, $v= 0 $ on $ \partial \Omega \setminus
D$, and  with $ \| \nabla v \| _ { L ^ 2 ( \Omega) }
\leq C \| f\|_ {L^2 ( \Omega)}$. This uses our assumption
(\ref{DOpen}) that $D$ is large.  Integrating by parts, writing
derivatives of $ \mu^\alpha_{ i}(u) $ in terms of $ \epsilon(u)$, and
integrating by parts again gives us 
\begin{equation}
\begin{split}
\int _ { \Omega } \mu^\alpha _ {i  } (u) f \, dy & =
\int _ { \Omega } \mu_i^{\alpha } (u) \div v \, dy \\
& = - \int _\Omega  \frac \partial { \partial y _ k} \mu^\alpha _ {i
 } (u) v_k \, dy  \\
& = - \int _ \Omega ( \frac \partial { \partial y _ i } \epsilon_k^
 {  \alpha } ( u ) - \frac \partial { \partial y _ \alpha }
\epsilon^i_ {k} (u) )v_k  \, dy \\
& = \int _ { \Omega } \epsilon_k^\alpha (u) \frac { \partial v _ k
}{ \partial y _ i } - \epsilon^i _{k}(u) \frac { \partial v _ k }{
  \partial y _ \alpha } \, dy .
\end{split}
\end{equation}
In the above calculation, the  boundary terms vanish because $ u$ is in $C^ \infty
_D(\bar  \Omega)$ and thus vanishes on a neighborhood of $D$, while $v$ 
vanishes on $ \partial \Omega \setminus D$.   If we take the supremum
over all $f$ with $ \| f \| _{L^ 2( \Omega)} \leq 1$, we obtain
\[
\| \mu^ \alpha _ { i  } (u) \| _{L^ 2( \Omega)} \leq C \| \epsilon (u)
\| _ {L^ 2( \Omega)}.
\]
As the constant in Proposition \ref{Bogo} may taken to be depend only
on 
$M$, the same holds for the constant in the Korn inequality of
this proposition. 
\end{proof}

\note{ 
Variable coefficient equations. 

Pressure term.

Dependence of constants. mostly fixed. 

Definition of starshaped domains.

To do: Variable coefficient operators.

Representation for $p$. 

Too many r's.  Mostly fixed. 

Give proof of symmetry. 

Estimates for $\delta^ \gamma G$  as in Mitrea-squared.  Omit or give
proof from TKB. 

Use of summation convention and definition of norm of a matrix.
}


\def\cprime{$'$} \def\cprime{$'$} \def\cprime{$'$} \def\cprime{$'$}
  \def\cprime{$'$}

\note{ Fix italics in Comm in PDE in reference to Taylor, Kim and Brown }

\end{document}